\documentclass[12pt]{article}
\usepackage{amsmath, amssymb, amsthm, delarray, epsfig, longtable, pstricks,
pst-plot, subfigure, tabularx, xspace} 

\theoremstyle{plain}
\newtheorem*{theorem*}{Theorem}
\newtheorem{theorem}{Theorem}

\begin{document}

\author{Abdulrahman Abdulaziz}
\title{The Chaos Game Versus Uniform Rotation: From Sierpinski Gaskets to
Periodic Orbits}
\maketitle

\abstract{ In this paper, we introduce a couple of dynamical systems that are
related to the Chaos Game. We begin by discussing different methods of
generating the Sierpinski gasket. Then we show how the transition from random to
uniform selection reduces the Sierpinski gasket to simple periodic orbits.
Next, we provide a simple formula for the attractor of each of the introduced
dynamical systems based only on the contraction ratio and the regular $n$-gon on
which the game is played. Finally, we show how the basins of attraction of a
particular dynamical system can generate some novel motifs that can tile the
plane.}

\section{The Chaos Game} 
\begin{figure}[ht]
\centering 
\includegraphics[width=6cm,height=4.85175cm]{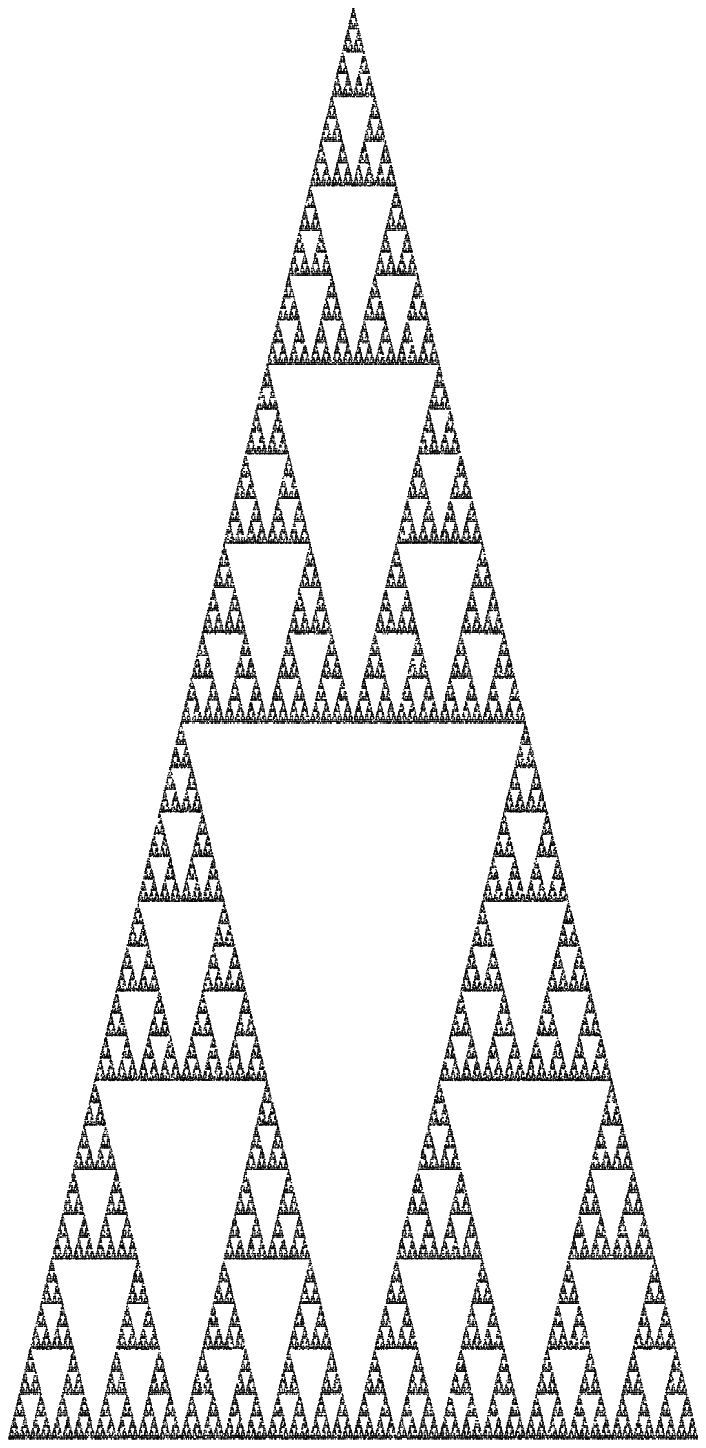}
\caption{The Sierpinski gasket.} \label{Fig:Srp}
\end{figure}

Let $\mathbf T=ABC$ be an equilateral triangle in the
plane. To play the Chaos Game (CG), start with any point in the plane, go
halfway from that point to a randomly chosen vertex of $\mathbf T$, and then
repeat the process, see~\cite{Barnsley:Fractal,Peitgen:Chaos}. Skipping the first thousand or
so iterates then plotting the points produced by the game, we obtain the
Sierpinski triangle or gasket shown in Figure~\ref{Fig:Srp}. The gasket, which
is called the \emph{attractor} of the game, was originally constructed by Waclaw
Sierpinski~\cite{Sierpinski} by removing the middle third of the original
triangle and then repeatedly removing the middle third of the remaining
triangles.

The Sierpinski triangle can be generated in yet another way described in
\cite{Schroeder:Fractal}. Start with a point $P$ in $\mathbf T$ and double the
distance from $P$ to the nearest vertex of $\mathbf T$ along the line from $P$
to that vertex. Clearly, most points will leave $\mathbf T$ under this game, but
it just happens that points that stay trapped inside $\mathbf T$ are exactly
those that belong to the Sierpinski gasket. This, however, is not a practical
way to generate the gasket since the probability that a randomly picked point
will stay in $\mathbf T$ under this game is virtually zero.

\section{New Games}
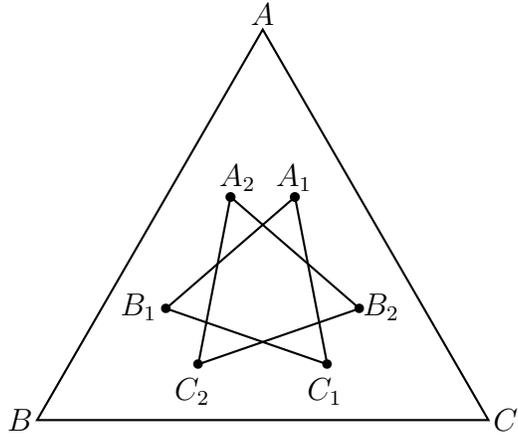
\begin{figure}[ht!]
\centering
\begin{pspicture}(-4.5,-0.25)(4.5,5.5)
\psset{xunit=0.75cm}
\psset{yunit=0.75cm}
\pspolygon(0,6.9282)(-4,0)(4,0)(0,6.9282)  
\pspolygon[showpoints=true](0.57143,3.959)(-1.7143,1.9795)(1.1429,0.98974)
(0.57143,3.959)
\pspolygon[showpoints=true](-0.57143,3.959)(-1.1429,0.98974)(1.7143,1.9795)
(-0.57143,3.959)
\rput(0,7.2){$A$}
\rput(-4.3, 0){$B$}
\rput(4.3, 0){$C$}
\rput(.55,4.3){$A_1$}
\rput(-2.2,2){$B_1$}
\rput(1.1,.5){$C_1$}
\rput(-.45,4.3){$A_2$}
\rput(2.1,2){$B_2$}
\rput(-1.25,.5){$C_2$}
\end{pspicture}
\caption{The attractors of the FVG.} \label{Fig:Str}
\end{figure}

Now suppose that instead of doubling the distance away from the nearest vertex,
we move halfway between the current point $P$ and that vertex of $\mathbf T$
which is farthest from $P$. Then, irrespective of where we start in the plane,
the attractor will either be the vertices of the triangle
$\mathbf T_r = A_1B_1C_1$ or those of the triangle $\mathbf T_l = A_2B_2C_2$,
see Figure~\ref{Fig:Str}.
We will show that the attractor reached depends solely on the location of the
initial point in the plane and consequently where it first lands in $\mathbf T$
as the game is played.


It turned out the attractors $\mathbf T_r$ and $\mathbf T_l$ can be obtained by
playing a uniform rotation game (URG). That is, instead of moving toward the
farthest vertex, we simply go around $\mathbf T$ in a clockwise or
counterclockwise direction. For counterclockwise rotation, first go halfway
between the current point and the vertex $A$, followed by $B$, and then $C$. In
this case, the attractor will be the vertices of the triangle $\mathbf T_r$. On
the other hand, clockwise rotation leads to the vertices of $\mathbf T_l$.

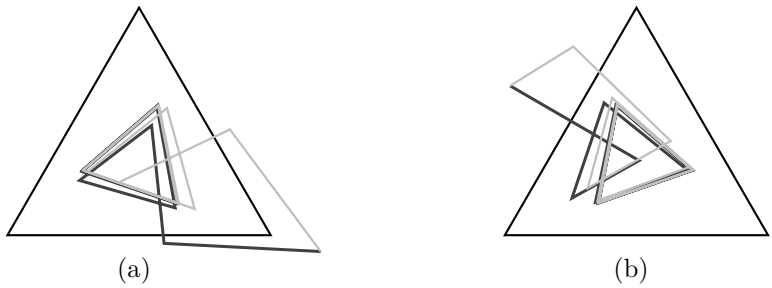
\begin{figure}[ht]
\centering
\psset{xunit=3.5cm}
\psset{yunit=3.5cm}
\subfigure[]{
\begin{pspicture}(-0.6,-0.3)(0.6,0.6)
\pspolygon(0.0, 0.577) (0.5, -0.288) (-0.5, -0.288)
\psline[linecolor=darkgray,linewidth=0.04]
(0.69,-0.35) (0.095,-0.319) (0.0475,0.129) (-0.226,-0.0799) (0.137,-0.184)
(0.0685,0.197) (-0.216,-0.0461) (0.142,-0.167) (0.0711,0.205)
(-0.215,-0.0419) (0.143,-0.165) (0.0714,0.206) (-0.214,-0.0413)
\psline[linecolor=lightgray,linewidth=0.03]
(0.69,-0.35) (0.345,0.114) (-0.0775,-0.0875) (0.211,-0.188) (0.106,0.195)
(-0.197,-0.047) (0.151,-0.168) (0.0757,0.205) (-0.212,-0.042)
(0.144,-0.165) (0.072,0.206) (-0.214,-0.0413) 
\end{pspicture}
}
\hskip 2cm
\subfigure[]{
\begin{pspicture}(-0.6,-0.3)(0.6,0.6)
\pspolygon(0.0, 0.577) (0.5, -0.288) (-0.5, -0.288)
\psline[linecolor=darkgray,linewidth=0.04]
(-0.48,0.28) (0.01,-0.00434) (-0.245,-0.147) (-0.123,0.215) (0.189,-0.0366)
(-0.156,-0.163) (-0.0779,0.207) (0.211,-0.0407) (-0.145,-0.165)
(-0.0723,0.206) (0.214,-0.0412) (-0.143,-0.165) (-0.0716,0.206)
(0.214,-0.0413)
\psline[linecolor=lightgray,linewidth=0.03]
(-0.48,0.28) (-0.24,0.429) (0.13,0.07) (-0.185,-0.109) (-0.0925,0.234)
(0.204,-0.0273) (-0.148,-0.158) (-0.0741,0.21) (0.213,-0.0395)
(-0.144,-0.164) (-0.0718,0.207) (0.214,-0.041) (-0.143,-0.165)
(-0.0715,0.206) (0.214,-0.0412)
\end{pspicture}
}
\caption{(a) Counterclockwise rotation (b) clockwise rotation.}
\label{Fig:Cun1}
\end{figure}

It is no accident that attractors of the FVG coincide with those of the URG. But
before we give the reason behind this, it should be noted that in the FVG, the
attractor could be $\mathbf T_r$ or $\mathbf T_l$, depending on the starting
point. This is in contrast with uniform rotation where the attractor depends
only on the choice of rotation. Now to illustrate how the two dynamical systems
lead to the same attractor, let us take a point that is attracted to
$\mathbf T_r$ under the FVG and by counterclockwise rotation. If we follow the
orbit of such a point we see that the two systems behave strikingly different
outside $\mathbf T$, but once inside $\mathbf T$ they both converge to
$\mathbf T_r$ quite fast, see Figure~\ref{Fig:Cun1}(a). Obviously, dark lines
represent the FVG while faded lines represent the URG. A similar phenomenon
occurs if we take a point that is attracted to $\mathbf T_l$ under clockwise
rotation, see Figure~\ref{Fig:Cun1}(b).

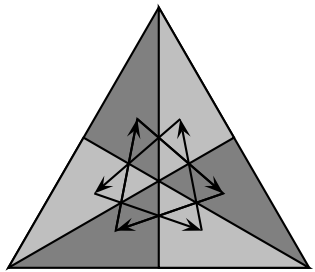
\begin{figure}[ht]
\centering
\psset{xunit=0.5cm}
\psset{yunit=0.5cm}
\begin{pspicture}(-4,0)(4,7)
\pspolygon[fillstyle=solid, fillcolor=gray](-4,0)(4,0)(0,6.9282)(-4,0)   
\pspolygon[fillstyle=solid, fillcolor=lightgray](-4,0)(0,2.3094)(-2,3.4641)(-4,0)
\pspolygon[fillstyle=solid, fillcolor=lightgray](0,0)(0,2.3094)(4,0)(0,0)
\pspolygon[fillstyle=solid,
fillcolor=lightgray](0,6.9282)(0,2.3094)(2,3.4641)(0,6.9282)
\psline[arrows=->, arrowsize=2pt 4](0.57143,3.959)(-1.7143,1.9795)
\psline[arrows=->, arrowsize=2pt 4](-1.7143,1.9795)(1.1429,0.98974)
\psline[arrows=->, arrowsize=2pt 4](1.1429,0.98974)(0.57143,3.959)
\pspolygon(-0.57143,3.959)(-1.1429,0.98974)(1.7143,1.9795)(-0.57143,3.959)
\psline[arrows=->, arrowsize=2pt 4](-1.1429,0.98974)(-0.57143,3.959)
\psline[arrows=->, arrowsize=2pt 4](-0.57143,3.959)(1.7143,1.9795)
\psline[arrows=->, arrowsize=2pt 4](1.7143,1.9795)(-1.1429,0.98974)
\end{pspicture}
\caption{Eventual behavior of the FVG.} \label{Fig:Dyn}
\end{figure}

Indeed, taking a closer look at this, we find that the two dynamical systems are
identical inside $\mathbf T$. To see this, let us divide $\mathbf T$ into dark and
light regions as in Figure~\ref{Fig:Dyn}. Then once inside $\mathbf T$ (not on one
of the bisectors), it is clear that the FVG is nothing but a clockwise rotation
for points in the dark area and a counterclockwise rotation for points in the
light area. If it happens that the current point is equidistant from two
vertices, then we break the tie arbitrarily and return to the previous case.

\section{Periodic Points}
In this section, we shall study periodic points of the URG, which are in turn
periodic points of the FVG. This is so since all points in the plane end up
inside $\mathbf T$ under both games and since for points in $\mathbf T$ the two
games are indistinguishable.

First, we prove that all points in the plane are attracted to the three vertices
of $\mathbf T_r$ under counterclockwise rotation. To do so, let us assume (without
loss of generality) that $\mathbf T$ is centered at the origin with $AB=1$.  Then
the vertices of $\mathbf T$ will have the coordinates
\[
A \left( 0, \frac{\sqrt{3}}{3} \right), \ B \left( -\frac{1}{2},
-\frac{\sqrt{3}}{6}\right), \ C \left( \frac{1}{2},
-\frac{\sqrt{3}}{6} \right),
\]
and the vertices of $\mathbf T_r$ and $\mathbf T_l$ will be\footnote{The vertices of
$\mathbf T_r$ and $\mathbf T_l$ are not written in simplest in order to easily
identify the similarities between their coordinates.}
\[
A_1 \left(\frac{1}{14} \, \raisebox{.7mm}{,} \ 5\frac{\sqrt{3}}{42}
\right), \ B_1 \left( -\frac{3}{14}, -\frac{\sqrt{3}}{42} \right), \ C_1
\left( \frac{2}{14}, -4\frac{\sqrt{3}}{42} \right),
\]
\[
A_2 \left(-\frac{1}{14}, 5\frac{\sqrt{3}}{42} \right), \ B_2 \left(
\frac{3}{14}, -\frac{\sqrt{3}}{42} \right), \ C_2 \left( -\frac{2}{14},
-4\frac{\sqrt{3}}{42} \right).
\]
It can be easily checked that
\[
A_1\overset{B}\longrightarrow B_1 \overset{C}\longrightarrow C_1
\overset{A}\longrightarrow A_1 \quad \text{and} \quad
A_2\overset{C}\longrightarrow B_2 \overset{B}\longrightarrow C_2
\overset{A}\longrightarrow A_2,
\]
where $A_1\overset{B}\longrightarrow B_1$ means $A_1$ goes to $B_1$ by moving
halfway toward $B$. In fact, it can easily checked that if
$P_0(x_0,y_0) \overset{A}\longrightarrow P_1$ and
$P_1\overset{B}\longrightarrow P_2$, then
\[
P_1\left(\frac{x_0}{2},\frac{y_0}{2}+\frac{\sqrt{3}}{6}\right) \quad \text{and}
\quad P_2\left(\frac{x_0}{4}-\frac{1}{4},\frac{y_0}{4}\right).
\]
It follows that 
\[
|P_1A_1|^2 - |P_2B_1|^2 = \frac{3}{16}
\left[\left(x_0-\frac{1}{7}\right)^2+\left(y_0+\frac{2}{7
\sqrt{3}}\right)^2\right],
\]
which is always positive except when $P_0=C_1$. Similarly, if
$P_2\overset{C}\longrightarrow P_3$, then
\[
|P_2B_1|^2 - |P_3C_1|^2 > 0
\]
unless $P_0=A_1$.  In other words, once we start rotating inside $\mathbf T$,
the distance between the current point and the nearest vertex of $\mathbf T_r$
will always be greater than the distance between its image and the nearest
vertex of $\mathbf T_r$. As we keep iterating, the system will tend to
$\mathbf T_r$, which makes it the sole attractor of the URG under
counterclockwise rotation. Similarly, it can be shown that $\mathbf T_l$ is the
attractor of the game under clockwise rotation.

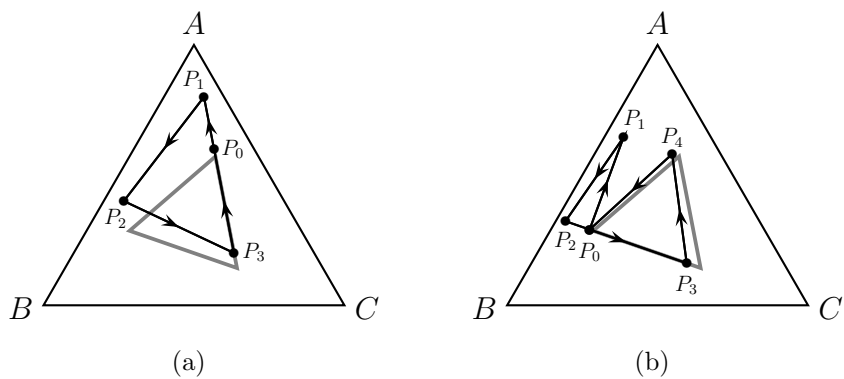
\begin{figure}[ht!]
\centering
\subfigure[]{
\begin{pspicture}(-2.5,-1.5)(2.5,3)
\psset{xunit=4cm}
\psset{yunit=4cm}
\psset{arrowscale=1.5,ArrowFill=true}
\pspolygon(0,0.577)(-0.500,-0.289)(0.500,-0.289)
\pspolygon[linewidth=0.05,linecolor=gray](0.0714,0.206)(-0.214,-0.0412)(0.143,-0.165)
\pspolygon[showpoints=true](0.0667,0.231)(0.0333,0.404)(-0.233,0.0577)(0.133,-0.115)
\psline[ArrowInside=->,
ArrowInsidePos=0.5](0.0667,0.231)(0.0333,0.404)(-0.233,0.0577)
(0.133,-0.115)(0.0667,0.231)
\rput(0,0.65){$A$}
\rput(-0.575, -0.3){$B$}
\rput(0.575, -0.3){$C$}
\rput(0.13,0.231){\scalebox{0.7}{$P_0$}}
\rput(0,0.46){\scalebox{0.7}{$P_1$}}
\rput(-0.26,0){\scalebox{0.7}{$P_2$}}
\rput(0.2,-0.11){\scalebox{0.7}{$P_3$}}
\end{pspicture}
}
\hskip 0.75cm
\subfigure[]{
\begin{pspicture}(-2.5,-1.5)(2.5,3)
\psset{xunit=4cm}
\psset{yunit=4cm}
\psset{arrowscale=1.5,ArrowFill=true}
\pspolygon(0,0.577)(-0.500,-0.289)(0.500,-0.289)
\pspolygon[linewidth=0.05,linecolor=gray](0.0714,0.206)(-0.214,-0.0412)(0.143,-0.165)
\pspolygon[showpoints=true](-0.226,-0.0372)(-0.113,0.270)(-0.306,-0.00931)(0.0968,-0.149)
(0.0484,0.214)(-0.226,-0.0372)
\psline[ArrowInside=->,
ArrowInsidePos=0.5](-0.226,-0.0372)(-0.113,0.270)(-0.306,-0.00931)(0.0968,-0.149)
(0.0484,0.214)(-0.226,-0.0372)
\rput(0,0.65){$A$}
\rput(-0.575, -0.3){$B$}
\rput(0.575, -0.3){$C$}
\rput(-0.226,-0.11){\scalebox{0.7}{$P_0$}}
\rput(-0.075,.33){\scalebox{0.7}{$P_1$}}
\rput(-0.306,-0.075){\scalebox{0.7}{$P_2$}}
\rput(0.0968,-0.22){\scalebox{0.7}{$P_3$}}
\rput(0.0484,0.27){\scalebox{0.7}{$P_4$}}
\end{pspicture}
}
\caption{A periodic point of (a) period 4 and (b) period 5.} \label{Uni:Per}
\end{figure}

Even though all points in the plane are attracted to $\mathbf T_r$ under
counterclockwise rotation, it is still possible to obtain periodic orbits of
periods other than three. Clearly, $A$ is the only fixed point of the game,
assuming we start the counterclockwise rotation at $A$. Also, $P_0(-1/3,0)$ is a
period 2 point since
\[
P_0\left(-\frac13,0\right) \overset{A}\longrightarrow P_1
\left(-\frac16,\frac{\sqrt{3}}{6}\right) \overset{B}\longrightarrow
P_0\left(-\frac13,0\right),
\]
where both $P_0$ and $P_1$ belong to the line $AB$. For a period $p>2$, the
periodic points fall inside $\mathbf T$.  For example, suppose that we want to
find a point of period four. Starting with the initial point $P_0(x_0,y_0)$ then
going counterclockwise, we get the following sequence of points
\begin{align*}
P_0(x_0,y_0)&\overset{A}\longrightarrow P_1\left(
\frac{x_0}{2}\raisebox{.1cm}{,} \, \frac{y_0}{2} + \frac{\sqrt{3}}{6} \right)
\overset{B}\longrightarrow P_2\left( \frac{x_0-1}{4} \raisebox{.1cm}{,}
\frac{y_0}{4} \right) \\
&\overset{C}\longrightarrow P_3\left(\frac{x_0+1}{8} \raisebox{.1cm}{,} \,
\frac{y_0}{8} - \frac{\sqrt{3}}{12} \right)\overset{A}\longrightarrow P_4\left(
\frac{x_0+1}{16} \raisebox{.1cm}{,} \, \frac{y_0}{16} + \frac{\sqrt{3}}{8}
\right).
\end{align*}
Hence, the desired point can be found by solving the linear system
$P_0=P_4$. This yields $x_0=1/15$, $y_0=2\sqrt{3}/15$, see Figure
\ref{Uni:Per}(a).
Observe that the nearest point on the orbit to $A_1$ is $P_0$, which lies on the
line $P_1P_3$ since
\[
P_3\overset{A}\longrightarrow P_0\overset{A}\longrightarrow P_1
\]
This will be true whenever the period is one more than a multiple of 3. In
Figure \ref{Uni:Per}(b), we show a point of period five. In this case, the point
$P_0$ must lie near $B_1$ on the line $B_1C$. This will always be true for
periods that are two more than a multiple of 3. However, there are no points
with a prime period that is a multiple of three. In this case, the point $P_0$
coincide with $C_1$. Otherwise, the image of $P_0$ will be closer to $A_1$ than
$P_0$ is to $C_1$, which we know is impossible.  Of course, we could obtain two
more periodic orbits of period $p$ if we move first toward $B$ or $C$. This is
equivalent to rotating the orbit obtained by moving first toward $A$ by an angle
of $2\pi/3$ or $4\pi/3$.

We can prove the above statements more rigorously by finding a formula for $P_n$
in general. First, observe that if
\[
M=
\begin{array}({rr})
\frac{1}{2}& 0 \\[3pt]
0& \frac{1}{2}
\end{array},
\
T_A=
\begin{array}({c})
0\\[3pt]
\frac{1}{2\sqrt{3}}
\end{array},
\
T_B=
\begin{array}({c})
\frac{-1}{4}\\[3pt]
\frac{-1}{4 \sqrt{3}}
\end{array},
\
T_C=
\begin{array}({c})
\ \ \frac{1}{4}\\[3pt]
\frac{-1}{4 \sqrt{3}}
\end{array}
\]
then going halfway from a point $P(x,y)$ to $A$, $B$ or $C$ amounts to adding
$T_A$, $T_B$ or $T_C$ to
\[
MP=
\begin{array}({cc})
\frac{1}{2}& 0 \\
0& \frac{1}{2}
\end{array}
\begin{array}({c}) x\\
y
\end{array}.
\]
It follows that for $P_0(x_0,y_0)$
\[
P_3 = M^3P_0 + M^2T_A+MT_B+T_C.
\]
Let $T=M^2T_A+MT_B+T_C$. Then continuing this way we get
\begin{equation}
P_{3k}=M^{3k}P_0 + \left(\sum _{l=0}^{k-1} M^{3 l}\right) T \label{Eqn:P3k}
\end{equation}
where $M^0$ is the identity matrix $I$. If $N=\left(I-M^{3k}\right)^{-1}$ and
$S=\sum _{l=0}^{k-1} M^{3 l}$, then solving $P_{3k}=P_0$ yields
\begin{equation}
P_0 =  NST
\end{equation}
It happens that in our case
\[
N = 
\begin{array}({cc})
\frac{2^{3k}}{2^{3k}-1}& 0 \\
0& \frac{2^{3k}}{2^{3k}-1} \\
\end{array},
\quad 
S=
\frac{8}{7}
\begin{array}({cc})
1-2^{-3k} & 0 \\
0 & 1-2^{-3k} \\
\end{array}
\]
and
\[
T=
\begin{array}({cc})
\ \ \frac{1}{8} \\[3pt]
-\frac{1}{4 \sqrt{3}}  \\
\end{array}.
\]
Therefore,
\[
NS =
\begin{array}({cc})
\frac{8}{7} & 0 \\
0 & \frac{8}{7} \\
\end{array}
\quad \text{and} \quad NST = \left(\frac{1}{7}, \, -\frac{2}{7 \sqrt{3}} \right)
=C_1.
\]
This proves that $C_1$ is the only periodic point of period $3k$ for
$k \in \mathbf N$. If $p=3k+1$, then solving $P_{3k+1}=P_0$, we obtain the
unique point
\[
\left(I-M^p\right)^{-1}(MST+T_A) = \left(\frac{2^{p-1}-1}{7
\left(2^p-1\right)},\, \frac{5\cdot 2^{p-1}+2}{7 \sqrt{3}
\left(2^p-1\right)}\right).
\]
It can be easily checked that the point lies on the line $A_1C_1$. Similarly,
for $p=3k+2$, we get the periodic point
\[
\left(I-M^p\right)^{-1}(M^2ST+MT_A+T_B) = \left(\frac{-3\ 2^{p-1}-1}{7
\left(2^p-1\right)},\, \frac{2-2^{p-1}}{7 \sqrt{3} \left(2^p-1\right)}\right),
\]
which falls on the line $B_1C_1$.

It is important to note that all periodic orbits, except the period 3 attractor,
are unstable and therefore they do not appear under iteration. In fact, they are
not even periodic points in the true sense of the word. For instance, if we
continue iterating the point of period 4 found earlier, we get
\[
P_0\overset{A}\longrightarrow P_1\overset{B}\longrightarrow
P_3\overset{C}\longrightarrow P_4\overset{A}\longrightarrow P_5=P_0
\overset{B}\longrightarrow P_6,
\]
where $P_6$ is different from $P_1$. In other words, even if a point $P_0$
returns to itself after a certain number of iterates, it will then get off track
and spiral toward $\mathbf T_r$. The reason behind this behavior is that the
game does not only depend on the contraction ratio but also on the chosen
vertex. For $P_6$ to be the same as $P_1$, we must move from $P_5$ toward $A$
instead of $B$, but then we would have violated the order of rotation. We
conclude that if we play the URG on any point in the plane, including pseudo
periodic points, we will always end up rotating around $\mathbf T_r$ or
$\mathbf T_l$ jumping from one vertex to the next, depending on whether the
rotation is counterclockwise or clockwise.

\section{Uniform Rotation on Regular $n$-gon}
The eventual behavior of the URG and FVG is the same only if the $n$-gon on
which the games are played is a triangle. For $n>3$, the URG yields periodic
orbits that form regular $n$-gons, while the FVG leads to star polygons. In this
section, we will study the URG, leaving the FVG for the next section.

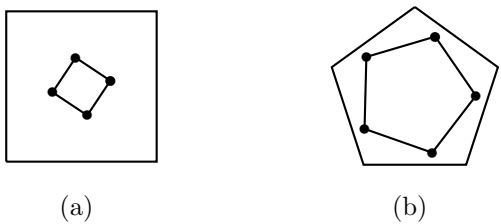
\begin{figure}[ht]
\centering
\psset{xunit=2cm}
\psset{yunit=2cm}
\subfigure[]{
\begin{pspicture}(-0.6,-0.6)(0.6,0.6)
\psline(-0.5,-0.5)(0.5,-0.5)(0.5,0.5)(-0.5,0.5)(-0.5,-0.5)
\psline[showpoints=true](0.192,0.0385)(-0.0385,0.192)(-0.192,-0.0385)
(0.0385,-0.192)(0.192,0.0385)
\end{pspicture}
}
\hskip 1.5cm
\psset{xunit=2.2cm}
\psset{yunit=2.2cm}
\subfigure[]{
\begin{pspicture}(-0.6,-0.5)(0.6,0.6)
\psline(0.,0.527)(-0.502,0.163)(-0.31,-0.427)(0.31,-0.427)(0.502,0.163)(0.,0.527)
\psline[showpoints=true](0.123,0.348)(-0.293,0.225)(-0.304,-0.21)(0.105,-0.354)
(0.369,-0.00945)(0.123,0.348)
\end{pspicture}
}
\caption{URG: (a) $n=4, \ r=1/3$ and (b) $n=5, \ r=2/3$.} \label{Uni:Game}
\end{figure}

First, observe that the inner triangles $\mathbf T_r$ and $\mathbf T_l$ are
tilted with respect to the original triangle $\mathbf T$ by an equal amount but
in a direction opposite to the direction of the game, see Figure \ref{Fig:Str}.
That is, if $\mathbf T_r$ is rotated clockwise by an angle $\alpha$, then
$\mathbf T_l$ is rotated counterclockwise by the same angle $\alpha$.  Moreover,
if the fraction of the distance traveled, henceforth called the
\emph{contraction ratio}, is changed from half the distance between the current
point and the chosen vertex to any other number between $0$ and $1$, then the
angle $\alpha$ changes as well.  This means that $\alpha$ does not only depend
on the number of vertices $n$, but also on the contraction ratio $r$, as shown
Figure \ref{Uni:Game}.

\begin{figure}[ht]
\centering
\psset{xunit=3cm}
\psset{yunit=3cm}
\begin{pspicture}(-0.7,-0.75)(0.81,1)
\psline(0.,0.851)(-0.809,0.263)(-0.5,-0.688)(0.5,-0.688)(0.809,0.263)(0.,0.851)
\psline[linestyle=dashed](0.21,0.51)(-0.42,0.357)(-0.469,-0.289)(0.13,-0.536)(0.55,
-0.0422)(0.21,0.51)
\psline[showpoints=true](0.55,-0.0422)(0.21,0.51)(0.,0.851)(0.809,0.263)
\psline[showpoints=true](0.,0.851)(0,0)(0.21,0.51)
\rput(0.62,0.05){$P$}
\rput(0.33,0.47){$Q$}
\rput(0,0.95){$V$}
\rput(0.925,0.263){$W$}
\rput(0.05,0.27){$\alpha$}
\rput(0.19,0.65){$\alpha$}
\rput(0,-0.1){$O$}
\end{pspicture}
\caption{The angle $\alpha$ for the URG on a regular
$n$-gon.} \label{Uni:AngFig}
\end{figure}
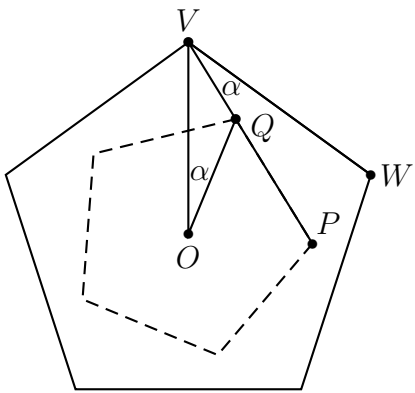

In order to calculate $\alpha_{n,r}$, we first find the inner $n$-gon whose
vertices are the attractor of the game. Suppose that $P(a,b)$ is a point on the
attractor and that $Q$ is the image of $P$ as we move toward a vertex $V$ of the
outer polygon $\mathbf P_n$, see Figure \ref{Uni:AngFig}. If $r$ is the fraction
of $PV$ traveled, then the point $Q$ is given by
\begin{equation}
Q=(1-r)P+rV
\end{equation}
Let $M_n$ be the $2\times2$ matrix that rotates a point in the plane by
$2\pi/n$. Then we can find $P$ by solving the equation
\begin{equation}
M_n^{n-1}Q = P. 
\end{equation}
In particular, if $V$ has the coordinates $(x,y)$ and $\rho=r-1$, then the
coordinates $(a,b)$ of $P$ are
\begin{equation}
(\rho+1)\left(\frac{\rho x + x\cos\frac{2\pi}{n} + y\sin\frac{2\pi}{n}}{1 +
\rho^2 + 2\rho\cos\frac{2\pi}{n}}, \ \frac{\rho y - x\sin\frac{2\pi}{n} +
y\cos\frac{2\pi}{n}}{1 + \rho^2 + 2 \rho\cos \frac{2\pi}{n}}
\right) \label{Uni:IniPt}
\end{equation}

Having found the coordinates of $P$ and thus $Q$, it is not difficult to see
that the angle by which the inner $n$-gon is rotated with respect to the outer
$n$-gon is simply $\alpha = \angle QOV$.
Moreover, since $W=M_n^{n-1}V$, $P=M_n^{n-1}Q$, and $P$, $Q$, and $V$ are
collinear, we must also have $\alpha = \angle PVW$.
It follows that if we use the values of $a$ and $b$ obtained in
(\ref{Uni:IniPt}), then both formulas for $\alpha$ should produce the same
answer. The tedious calculation\footnote{We can simplify the calculation of
$\alpha$ by assuming without loss of generality that the coordinates of $V(x,y)$
are $(0,1)$} yields
\begin{equation}
\alpha_{n,\rho} = \cos ^{-1}\left(\frac{1+\rho \cos \frac{2 \pi
}{n}}{\sqrt{1+\rho^2+2 \rho \cos \frac{2 \pi }{n}}}\right) \label{Uni:AngFor}
\end{equation}
Note that as $r\to1$, $\rho\to0$ and so $\alpha_{n,\rho}$ goes to
$\cos^{-1}(1)=0$. On the other hand, as $r\to0$, $\rho\to-1$ and we have
\begin{align*}
\lim_{\rho\to-1}\alpha_{n,\rho}& = \cos ^{-1}\left( \frac{2\sin^2
\frac{\pi}{n}}{\sqrt{4\sin^2 \frac{\pi}{n}}}\right) = \cos ^{-1}\left(\sin
\frac{\pi }{n}\right) \\[7pt]
&= \cos ^{-1}\left(\cos \left(\frac{\pi}{2} - \frac{\pi }{n}\right)\right) =
\frac{(n-2)\pi}{2n},
\end{align*}
which is half the angle of a regular $n$-gon.

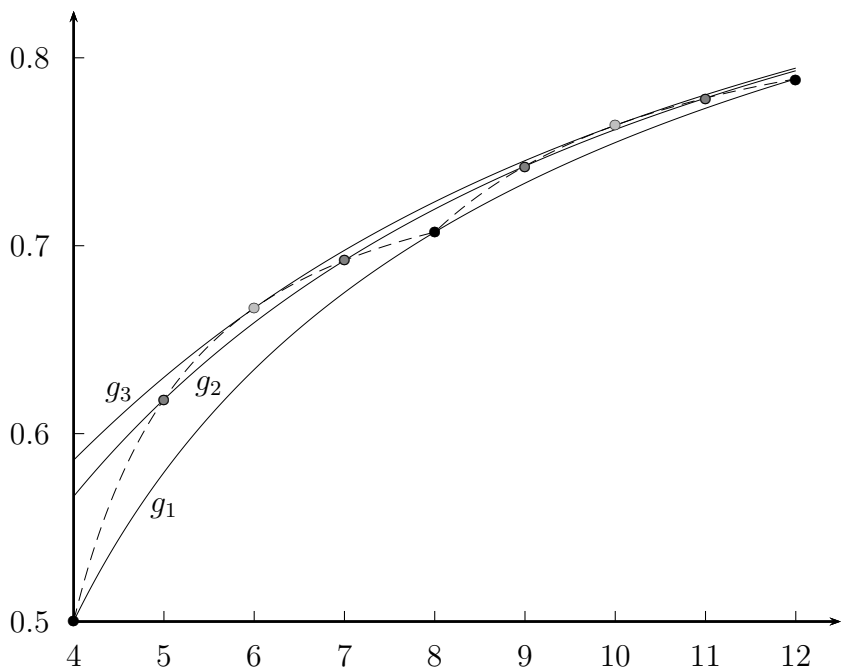
\begin{figure}[ht]
\psset{xunit=1.2cm}
\psset{yunit=25cm}
\centering
\savedata{\dataA}[{{4.,0.5}{8.,0.707}{12.,0.788}}]
\savedata{\dataB}[{{5.00,0.618}{9.00,0.742}{7.00,0.692}{11.0,0.778}}]
\savedata{\dataC}[{{6.00,0.667}{10.0,0.764}}]
\begin{pspicture}(3.5,0.48)(12.5,0.84)
\psaxes[tickstyle=top,ticks=all,labels=all,Ox=4,dy=0.1,Dy=0.1,Oy=0.5,Dx=1,dx=1]{->}(4,0.5)(4,0.5)(12.5,0.825)
\psplot[plotpoints=200,linewidth=0.01cm]{4}{12}{1 1 PI x div TAN add div}
\psplot[plotpoints=200,linewidth=0.01cm]{4}{12}{1 1 PI x div SIN add div}
\psplot[plotpoints=200,linewidth=0.01cm]{4}{12}{1 1 2 PI 2 x
mul div SIN mul add div} 
\psplot[plotpoints=200,linestyle=dashed,linewidth=0.01cm]{4}{12}{PI 2 PI x 4 div
floor mul mul add 
x div SIN PI x div SIN PI 2 PI x 4 div floor mul mul add x div SIN add div}
\listplot[linewidth=1pt,plotstyle=dots,dotstyle=o,fillcolor=black]{\dataA}
\listplot[linewidth=1pt,plotstyle=dots,dotstyle=Bo,fillcolor=gray]{\dataB}
\listplot[linewidth=1pt,plotstyle=dots,dotstyle=o,fillcolor=lightgray]{\dataC}
\rput(5,0.56){$g_1$}
\rput(5.5,0.625){$g_2$}
\rput(4.5,0.622){$g_3$}
\end{pspicture}
\caption{The intersection of $k_x$ with $g_1$, $g_2$ and
$g_3$.} \label{Fig:Kiss}
\end{figure}

Although the contraction ratio $r$ could be any number in $[0,1]$, it is more
natural to choose a ratio that produces a fractal similar to the Sierpinski
gasket under the CG. This particular ratio has been called the \emph{kissing
ratio}, see \cite{Schlicker,Bates}. In fact, it was shown in \cite{Aziz} that
the kissing ratio is
\[
r=k_n = \frac{\sin \frac{\pi+2 \pi \left\lfloor
\frac{n}{4}\right\rfloor}{n}}{\sin \frac{\pi+2 \pi \left\lfloor
\frac{n}{4}\right\rfloor}{n}+\sin \frac{\pi }{n}} =
\begin{cases}
\dfrac{1}{1+\tan \frac{\pi}{n}}& \text{if } n \equiv 0 \bmod 4 \\[15pt]
\dfrac{1}{1+2\sin \frac{\pi}{2n}}& \text{if } n \equiv 1,3 \bmod 4 \\[15pt]
\dfrac{1}{1+\sin \frac{\pi}{n}}& \text{if } n \equiv 2 \bmod 4
\end{cases}  
\]
where $n$ is the number of sides in the regular polygon. Moreover, if we define
\[
k_x = \frac{\sin \frac{\pi+2 \pi \left\lfloor \frac{x}{4}\right\rfloor}{x}}{\sin
\frac{\pi+2 \pi \left\lfloor \frac{x}{4}\right\rfloor}{x}+\sin \frac{\pi }{x}}
\]
for any positive real number $x$, then we can visually see how 
$k_x$ intersects with the functions
\[
g_1(x)=\dfrac{1}{1+\tan \frac{\pi}{x}}, \quad g_2(x)=\dfrac{1}{1+2\sin
\frac{\pi}{2x}}, \quad \text{and} \quad g_3(x)=\dfrac{1}{1+\sin \frac{\pi}{x}}
\]
at only integer values, as illustrated in Figure \ref{Fig:Kiss}.

\begin{figure}[ht]
\psset{xunit=0.5cm}
\psset{yunit=5cm}
\centering
\begin{pspicture}(1,-0.12)(14,0.65)
\psaxes[tickstyle=top,ticks=all,labels=all,dy=0.25,Dy=0.25,Ox=3,Dx=3,dx=3]{->}(3,0)(3,0)(14,0.6)
\psline[showpoints=true](3, 0.333473)(4, 0.463648)(5, 0.390713)(6, 0.333473)(7,
0.289632)(8, 0.255495)(9, 0.203602)(10, 0.169861)(11,0.14618)(12, 0.128618)(13,
0.108593)
\end{pspicture}
\caption{Tilt angle for the URG with $r=k_n$.} \label{Uni:RotAng}
\end{figure}
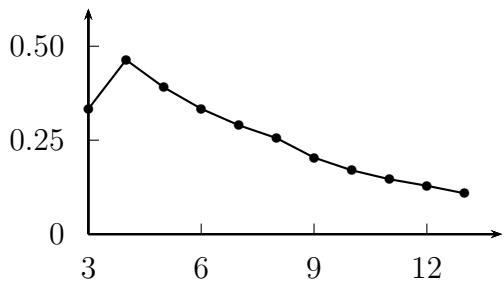

Clearly, for $x\geq3$ we have
\[
g_1(x)<g_2(x)<g_3(x).
\]
Note how $g_1$ touches $k_x$ only from below, $g_2$ crosses $k_x$, and $g_3$
touches $k_x$ only from above. If we set $\kappa_n=k_n-1$, then for a uniform
rotation with a contraction ratio $k_n$, we have
\begin{equation}
\alpha_n = \cos ^{-1} \left( \frac{1+\kappa_n \cos \frac{2 \pi
}{n}}{\sqrt{1+\kappa_n^2+2 \kappa_n \cos \frac{2 \pi }{n}}}
\right) \label{UniAng:Formula} 
\end{equation}
Here $\alpha_n$ is a function of $n$ only since the kissing ratio is a function
of $n$ as well. Curiously, the angle $\alpha_n$ increases as $n$ goes from $3$
to $4$ then decreases from $n=4$ on, tending to zero as $n$ goes to infinity,
see Figure~\ref{Uni:RotAng}. Another curiosity is that on the way down from the
maximum at $n=4$, when we reach $n=6$ we get $\alpha_6=\alpha_3$, which is the
only case where two distinct integers yield the same value of $\alpha$.

Next, we will show that the sum of $\alpha_{n,r}$ diverges for $0<r<1$.  Using
L'H\^opital's rule, it can be proved that
\[
\lim_{n\to\infty}\frac{\alpha_{n,r}}{\frac1n}= \frac{0}{0} = 2 \pi
\left(\frac{1}{r}-1\right),
\]
which is positive for $r\in (0,1)$. By the limit comparison test, it follows
that
\[
\sum_{n=3}^{\infty} \alpha_{n,r} = \infty.
\]
On the other hand, if we replace $r$ by the kissing ration $k_n$, then the sum
$\alpha_n$ converges since in this case $k_n\to1$ as $n\to \infty$. This can be
proved by splitting the series into three series depending on the remainder of
the division of $n$ by 4.  It turned out
\[
\lim_{n\to\infty}\frac{\alpha_{n}}{\frac{1}{n^2}} =
\begin{cases}
\frac{\pi^2}{8}& \text{if } n = 2k \\[5pt]
\frac{\pi^2}{2}& \text{if } n = 2k+1 
\end{cases}
\]
Hence, the sum of the whole series converges, again by the limit comparison
test.

\section{Farthest Vertex on Regular $n$-gon}
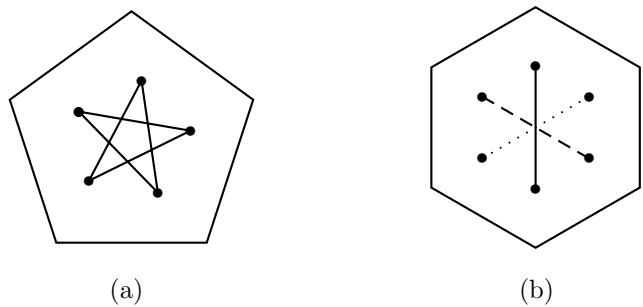
\begin{figure}[ht]
\centering
\psset{xunit=2cm}
\psset{yunit=2cm}
\subfigure[]{
\begin{pspicture}(-0.6,-0.8)(0.6,0.9)
\pspolygon(0.,0.851)(-0.809,0.263)(-0.5,-0.688)(0.5,-0.688)(0.809,0.263)
\pspolygon[showpoints=true](-0.35,0.184)(0.175,-0.355)(0.0669,0.39)
(-0.283,-0.276)(0.392,0.0569)(-0.35,0.184)(-0.35,0.184)
\end{pspicture}
}
\hskip 2cm
\psset{xunit=1.6cm}
\psset{yunit=1.6cm}
\subfigure[]{
\begin{pspicture}(-1.1,-1.1)(1.2,0.6)
\pspolygon(0.,1.)(-0.866,0.5)(-0.866,-0.5)(0.,-1.)(0.866,-0.5)(0.866,0.5)
\psline[showpoints=true](0.,0.513)(0.,-0.513)
\psline[showpoints=true,linestyle=dashed](-0.444,0.257)(0.444,-0.257)
\psline[showpoints=true,linestyle=dotted](-0.444,-0.257)(0.444,0.257)
\end{pspicture}
}
\caption{Attractors of the FVG (a) $n$ odd and (b) $n$ even.}
\label{Far:Stars}
\end{figure}

For $n > 3$, the FVG yields different attractors than the URG. In fact, the
eventual behavior of the FVG depends on whether $n$ is even or odd. More
precisely, for $n$ odd, the attractors are still periodic orbits of $n$ points
each; but instead of drawing a regular polygon, the line segments joining
consecutive vertices of an attractor now draw a star polygon, see Figure
\ref{Far:Stars}(a). On the other hand, for $n$ even, we get $n/2$ attractors,
each consisting of two points that lie on a line joining opposite vertices of
the original $n$-gon, as shown in Figure \ref{Far:Stars}(b).

Suppose that $n$ is odd, $P$ is a vertex of the attracting star polygon, and $Q$
is the image of $P$ as we travel toward the vertex $V(x,y)$ of the outer polygon
$\mathbf P_n$, see Figure \ref{Star:Attrc}. First, note that if $n$ is odd and
$V$ is chosen on the $y$-axis, then the two vertices that are farthest from $V$
lie on the horizontal side $WW'$ of $\mathbf P_n$. In addition, if $m=(n-1)/2$,
then the $W$ and $W'$ are exactly $m$ vertices from $V$, depending on whether we
move clockwise or counterclockwise. This means that moving counterclockwise from
$V$, we reach $W'$ after a rotation by $\pi-\pi/n$ and $W$ after a rotation by
$\pi+\pi/n$. It follows that we can get back from $Q$ to $P$ by a rotation $M$
of $\pi+\pi/n$, and so the coordinates of $P$ can be found by solving the
equation $P=MQ$. In this case, $P$ is given by
\begin{equation}
(\rho+1)\left( \frac{\rho x - x\cos\frac{\pi}{n} + y\sin\frac{\pi }{n}} {1 +
\rho ^2 - 2\rho\cos\frac{\pi}{n}},\ \frac{\rho y - x \sin \frac{\pi }{n} - y\cos
\frac{\pi }{n}}{1+\rho^2-2 \rho \cos \frac{\pi }{n}} \right) \label{Star:IniPt}
\end{equation}

As for the URG, the angle by which star polygon is rotated with respect to the
outer $n$-gon is $\alpha = \angle QOV = \angle PVW$. But $W$ in this case is the
farthest clockwise vertex from $V$, while it was the closest clockwise vertex to
$V$ in the URG. Knowing that $W=MV$, we can calculate either $\angle QOV$ or
$\angle PVW$ to obtain
\begin{figure}[ht]
\centering
\psset{xunit=3cm}
\psset{yunit=3cm}
\begin{pspicture}(-0.7,-0.81)(0.81,1)
\psline(0.,0.851)(-0.809,0.263)(-0.5,-0.688)(0.5,-0.688)(0.809,0.263)(0.,0.851)
\psline[linestyle=dashed](-0.35,0.184)(0.175,-0.355)(0.0669,0.39)(-0.283,-0.276)
(0.392,0.0569)(-0.35,0.184)(-0.35,0.184)
(0.392,0.0569)(-0.35,0.184)
\psline[showpoints=true](0.175,-0.355)(0.0669,0.39)(0.,0.851)(0.5,-0.688)
\rput(0.175,-0.45){\scalebox{0.8}{$P$}}
\rput(-0.025,0.4){\scalebox{0.8}{$Q$}}
\rput(0,0.93){\scalebox{0.8}{$V$}}
\rput(0.5,-0.8){\scalebox{0.8}{$W$}}
\rput(-0.5,-0.8){\scalebox{0.8}{$W'$}}
\rput(0.15,0.2){$\alpha$}
\rput(0,-0.08){\scalebox{0.7}{$O$}}
\psdots(0,0)(-0.5,-0.688)
\end{pspicture}
\caption{The angle $\alpha$ for the FVG on a regular
$n$-gon.} \label{Star:Attrc}
\end{figure}
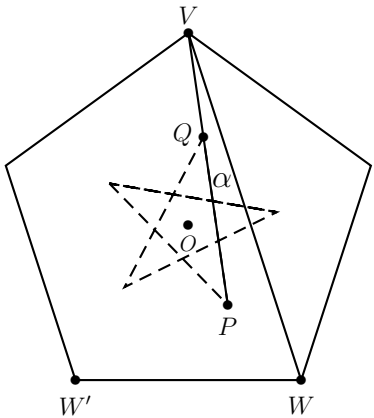
\begin{equation}
\alpha_{n,\rho}=\cos ^{-1}\left(\frac{1-\rho \cos \frac{\pi }{n}}{\sqrt{ 1 +
\rho ^2 -2\rho\cos\frac{\pi}{n}}}\right) \label{Star:Ang}
\end{equation}

Of course the other inner star will be rotated in a clockwise direction. As in
the case of the triangle, the choice of rotation depends on where inside the
$n$-gon the system lands first. In Figure \ref{Star:Attrc}, if we assume that
$P$ is the first iterate that falls inside the pentagon, then clearly the game
will continue in counterclockwise direction. On the other hand, the reflection
of $P$ with respect to the $y$-axis will lead to clockwise rotation. Moreover,
counterclockwise rotation will force the inner star to be rotated clockwise with
respect to the outer polygon and vice versa.

Having found the formulas for $P$ and $\alpha$ for the FVG when $n$ is odd, it
only is natural to compare them with those obtained for the URG. In fact, it is
noticeable that one can get from the URG to the FVG simply by replacing every
occurrence of $\cos 2\pi/n$ by $-\cos \pi/n$. From this we see that when $n=3$,
the two games yield the same attractors since $\cos 2\pi/3=-\cos\pi/3$.

\begin{figure}[ht]
\psset{xunit=0.5cm}
\psset{yunit=5cm}
\centering
\begin{pspicture}(1,-0.12)(14,0.5)
\psaxes[tickstyle=top,ticks=all,labels=all,dy=0.25,Dy=0.25,Ox=3,Dx=3,dx=3]{->}(3,0)(3,0)(14,0.45)
\psline[showpoints=true](3,0.333)(4,0.255)(5,0.170)(6,0.129)(7,0.104)(8,0.0880)(9,0.0709)
(10.0,0.0595)(11.0,0.0514)(12.0,0.0454)(13.0,0.0391)
\end{pspicture}
\caption{Tilt angle for the FVG with $r=k_n$.} \label{Star:RotAng}
\end{figure}
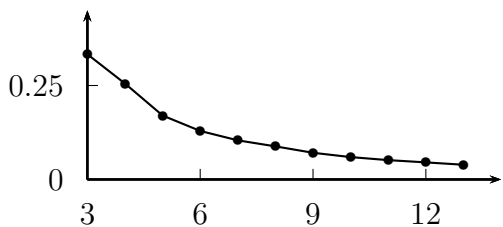

Another similarity between the two games is that when $r\to 1$, $\alpha$ tends
to $0$ for all $n$. Also, as in the URG, the sum of $\alpha_{n,r}$ diverges for
any fixed $r$, while the sum converges if we take $r=k_n$. However, as $r\to 0$,
we get
\[
\lim_{\rho\to-1}\alpha_{n,r} = \cos ^{-1}\left( \frac{2\cos^2
\frac{\pi}{2n}}{\sqrt{4\cos^2 \frac{\pi}{2n}}}\right) = \cos ^{-1}\left(\cos
\frac{\pi }{2n}\right) = \frac{\pi}{2n},
\]
which is different from the value of $\alpha_{n,r}$ in the URG. Moreover, in
this case, $\alpha_{n+1} < \alpha_n$ for $n\geq 3$, as shown in Figure
\ref{Star:RotAng}.

Finally, if $n$ is even, then we get $n/2$ attractors, each consisting of a pair
of points falling on a diagonal joining opposite vertices of the outer
$n$-gon. In other words, if $Q$ is the image of $P$ as we move toward $V(x,y)$,
then $P=MQ$ is the image of $Q$, where $M$ is rotation by $\pi$. In this case,
the coordinates of $P$ are
\begin{equation}
\frac{\rho + 1}{\rho -1}\left(x,y\right) \label{Star:Even}
\end{equation}
Note that the coordinates of $P$ are now independent of $n$. Also, the inner
polygon whose vertices are the $n/2$ attracting pairs of points is not rotated
with respect to the original $n$-gon, and so $\alpha=0$, see
Figure~\ref{Far:Stars}(b). In fact, we could have obtained the formulas for $P$
and $\alpha$ from the formulas for $n$ odd by replacing $\pi/n$ by 0.

\section{Transformation Matrix}
The counterclockwise rotation game on $\mathbf T$ with contraction ratio $1/2$
yielded the triangle $\mathbf T_r$, whose sides are $1/\sqrt{7}$ the sides of
$\mathbf T$. In general, if we take the regular $n$-gon to be centered at the
origin with $(x,y)$ as one of its vertices, then that the point $P$ of the
attractor that moves toward$(x,y)$ was given by (\ref{Uni:IniPt}).  Also, if we
calculate the length of the side of the inner polygon, which is the distance
between $P$ and its image under a rotation by $2\pi/n$, we get
\begin{equation}
l_{n,\rho}=\frac{2c(\rho+1)\sin\frac\pi n}{\sqrt{1 + \rho^2 + 2\rho\cos \frac{2 \pi
}{n}}}
\end{equation}
where $c=\sqrt{x^2+y^2}$ is the radius of the circumscribing circle.  Since the
length of the side of a regular polygon in a circle of radius $c$ is
$L_n=2c\sin \pi/n$, we get a scaling factor
\begin{equation}
\lambda_{n,\rho}=\frac{l_{n.\rho}}{L_n}=\frac{\rho+1}{\sqrt{1 + \rho^2 +
2\rho\cos \frac{2 \pi }{n}}}
\end{equation}

\begin{figure}[ht!]
\centering
\subfigure[]{
\begin{pspicture}(-2.5,-1.5)(2.5,3)
\psset{xunit=4cm}
\psset{yunit=4cm}
\psset{arrowscale=1.5,ArrowFill=true}
\pspolygon(0,0.577)(-0.500,-0.289)(0.500,-0.289)
\pspolygon(0.0405,0.0858)(-0.0946,-0.0078)(0.0541,-0.0780)(0.0405,0.0858)
\pspolygon(0.0714,0.206)(-0.214,-0.041)(0.143,-0.165)(0.0714,0.206)
\pspolygon(0.0769,0.311)(-0.308,-0.089)(0.231,-0.222)(0.0769,0.311)
\pspolygon(0.0526,0.456)(-0.421,-0.182)(0.368,-0.273)(0.0526,0.456)
\end{pspicture}
}
\hskip 0.75cm
\subfigure[]{
\begin{pspicture}(-2.5,-2.1)(2.5,2.5)
\psset{xunit=3.5cm}
\psset{yunit=3.5cm}
\psset{arrowscale=1.5,ArrowFill=true}
\pspolygon(0.5,0.5)(-0.5,0.5)(-0.5,-0.5)(0.5,-0.5)(0.5,0.5)
\pspolygon(0.140,0.020)(-0.020,0.140)(-0.140,-0.020)(0.020,-0.140)(0.140,0.020)
\pspolygon(0.300,0.100)(-0.100,0.300)(-0.300,-0.100)(0.100,-0.300)(0.300,0.100)
\pspolygon(0.400,0.200)(-0.200,0.400)(-0.400,-0.200)(0.200,-0.400)(0.400,0.200)
\pspolygon(0.480,0.360)(-0.360,0.480)(-0.480,-0.360)(0.360,-0.480)(0.480,0.360)
\end{pspicture}
}
\caption{The inner polygons for different ratios.} \label{Fig:InPoly}
\end{figure}
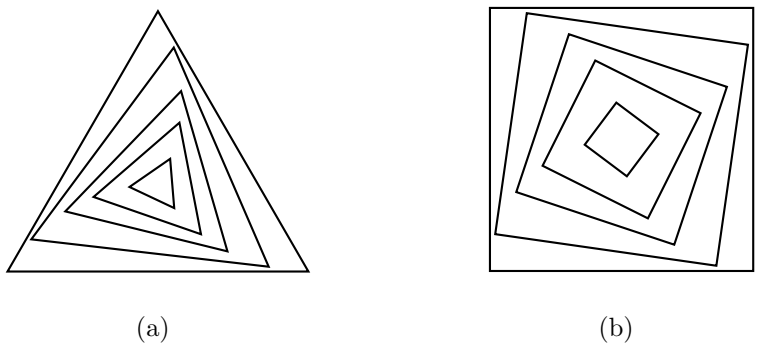

Now knowing the angle $\alpha_{n,\rho}$ and the scaling factor
$\lambda_{n,\rho}$, it follows that if
\[
M_{n,\rho}=
\begin{array}({rr})
\cos \alpha_{n,\rho}& -\sin \alpha_{n,\rho} \\
\sin \alpha_{n,\rho}&  \cos \alpha_{n,\rho}
\end{array}
\]
then $\lambda_{n,\rho}M_{n,\rho}$ transforms the original $n$-gon to the
attracting inner $n$-gon. The other inner $n$-gon is obtained y replacing
$\alpha_{n,\rho}$ by $-\alpha_{n,\rho}$. In Figure \ref{Fig:InPoly}, we show the
inner polygons for a triangle and a square for
\[
r = \frac14, \, \frac12\, \frac23, \, \frac67.
\]

Observe how the angle $\alpha_{n,\rho}$ decreases as $\lambda_{n,\rho}$
increases so that each extended side of the inner polygons passes exactly
through one vertex of the original polygon.  This extremely simple yet powerful
conclusion can be restated as an elegant
\begin{theorem}
Let $\mathbf P_n$ be a regular $n$-gon. Then the attractor of a counterclockwise
(clockwise) uniform rotation with contraction ratio $r$ is a periodic orbit of
period $n$ that can be obtained by applying the linear transformation
$\lambda M$ to the vertices of $\mathbf P_n$, where $M$ is a clockwise
(counterclockwise) rotation by
\[
\alpha_{n,\rho} = \cos ^{-1}\left(\frac{1+\rho \cos \frac{2 \pi
}{n}}{\sqrt{1+\rho^2+2 \rho \cos \frac{2 \pi }{n}}}\right),
\]
$\lambda$ is a scaling factor given by
\[
\lambda_{n,\rho}=\frac{1+\rho}{\sqrt{1 + \rho^2 + 2\rho\cos \frac{2 \pi }{n}}},
\]
and $\rho=r-1$.
\end{theorem}

\begin{figure}[ht!]
\centering
\psset{xunit=2cm}
\psset{yunit=2cm}
\subfigure[]{
\begin{pspicture}(-0.6,-0.8)(0.6,0.9)
\pspolygon(0.,0.851)(-0.809,0.263)(-0.5,-0.688)(0.5,-0.688)(0.809,0.263)
\pspolygon(-0.35,0.184)(0.175,-0.355)(0.0669,0.39)
(-0.283,-0.276)(0.392,0.0569)(-0.35,0.184)(-0.35,0.184)
\pspolygon[showpoints=true](0.0669,0.39)(-0.35,0.184)(-0.283,-0.276)(0.175,-0.355)
(0.392,0.0569)
\end{pspicture}
}
\hskip 2cm
\psset{xunit=1.6cm}
\psset{yunit=1.6cm}
\subfigure[]{
\begin{pspicture}(-1.1,-1.1)(1.2,0.6)
\pspolygon(0.,1.)(-0.866,0.5)(-0.866,-0.5)(0.,-1.)(0.866,-0.5)(0.866,0.5)
\psline(0.,0.513)(0.,-0.513)
\psline(-0.444,0.257)(0.444,-0.257)
\psline(-0.444,-0.257)(0.444,0.257)
\pspolygon[showpoints=true](0.,0.513)(-0.444,0.257)(-0.444,-0.257)(0.,-0.513)
(0.444,-0.257)(0.444,0.257)
\end{pspicture}
}
\caption{(a) Star polygon (b) rays from the origin.}
\label{Far:Polys}
\end{figure}
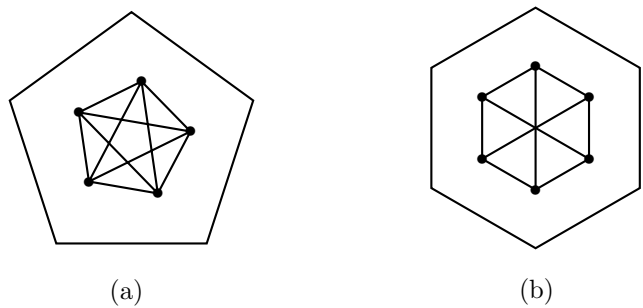

We can derive a similar theorem for the FVG by thinking of the attractors,
whether they are vertices of star polygons or endpoints of lines through the
origin, as diagonals inside an $n$-gon, see Figure~\ref{Far:Polys}.  If $n$ is
odd, then the scaling factor $\lambda$ can be obtained from the URG by replacing
$\cos(2\pi/n)$ by $-\cos(\pi/n)$, while if $n$ is even we replace $\cos(2\pi/n)$
by $-\cos(0)$. This yields

\begin{theorem}
Let $\mathbf P_n$ be a regular $n$-gon with $n$ odd. Then the attractor of the
FVG with contraction ratio $r$ is a periodic orbit of period $n$ that can be
obtained by applying the linear transformation $\lambda M$ to the vertices of
$\mathbf P_n$, where $M$ is a clockwise or counterclockwise rotation by
\[
\alpha_{n,\rho} = \cos ^{-1}\left(\frac{1-\rho \cos
\frac{\pi}{n}}{\sqrt{1+\rho^2-2 \rho \cos \frac{\pi}{n}}}\right),
\]
$\lambda$ is a scaling factor given by
\[
\lambda_{n,\rho}=\frac{1+\rho}{\sqrt{1 + \rho^2 - 2\rho\cos \frac{\pi}{n}}},
\]
and $\rho=r-1$.
\end{theorem}

Whether the rotation is clockwise or counterclockwise depends on the location of
the initial point in the plane and consequently the place where it first lands
in $\mathbf P_n$.

Finally, if $n$ is even, we get $n/2$ attractors, each consisting of a period 2
orbit. In this case, there is no rotation and the scaling factor is simply
\begin{equation}
\lambda_\rho = \frac{1 + \rho}{1 - \rho},
\end{equation}
which does not depend on $n$. Indeed, apart from the fact that the number of
attractors is now $n/2$ instead of $2$, the formulas for $\alpha$ and $\lambda$
we derived for odd $n$ yields the corresponding formula for even $n$, if $\pi/n$
in the formula for odd $n$ is replaced by zero.

\section{Basins of Attraction}
\begin{figure}[ht!]
\centering
\subfigure[]{
\includegraphics[width=3cm,height=3cm]{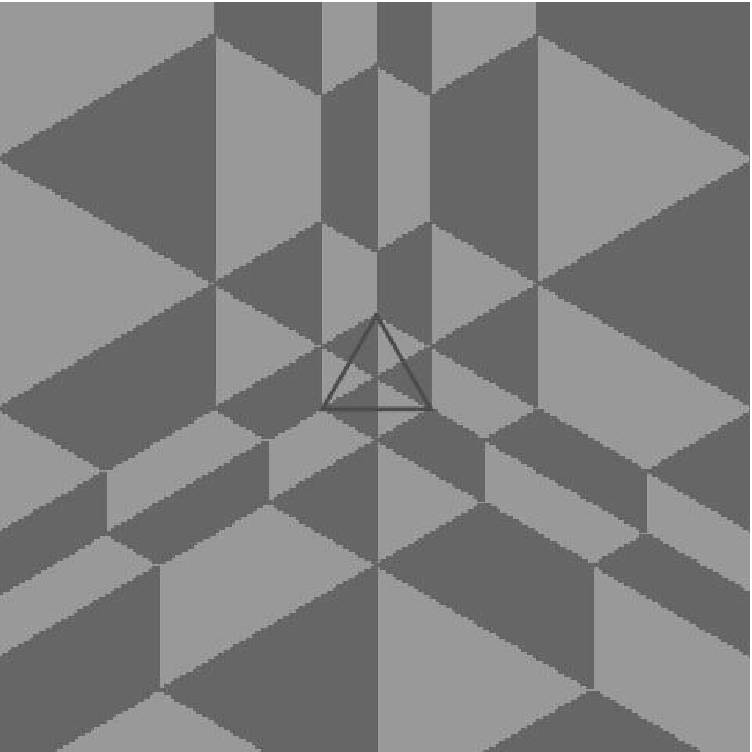} 
}
\hskip 2cm
\subfigure[]{
\includegraphics[width=3.05cm,height=3.05cm]{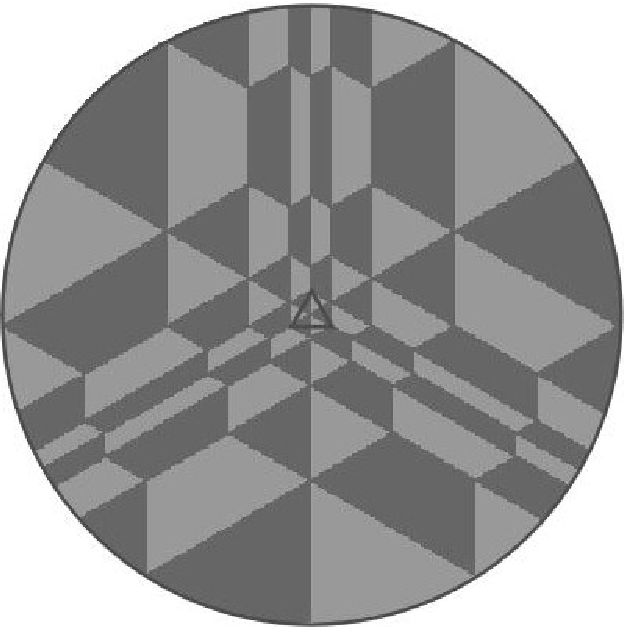} 
}
\caption{Basins for $n=3$ (a) square of side 7 (b) circle of radius
7.} \label{Basin:Odd3}
\end{figure}

For the URG, all points in the plane are attracted to $\mathbf T_r$ under
counterclockwise rotation and to $\mathbf T_l$ under clockwise rotation. On the
other hand, if we try to find the basins of attraction of $\mathbf T_r$ and
$\mathbf T_l$ under the FVG, we get the carpet shown in
Figure~\ref{Basin:Odd3}(a) with the triangle $\mathbf T$ in the
middle\footnote{Unless otherwise stated, the contraction ratio used is always
the kissing ratio.}. It is not difficult to see that light regions are attracted
to $\mathbf T_r$ while dark regions are attracted to $\mathbf T_l$. More
precisely, a point in a light (dark) region outside $\mathbf T$ will eventually
reach a light (dark) region inside $\mathbf T$, and thus spiral counterclockwise
(clockwise) toward $\mathbf T_r$ ($\mathbf T_l$). Like the original triangle
$\mathbf T$, the basins of attraction are invariant under rotation by $2k\pi/3$
for $k=1,2,3$, which can be clearly seen in Figure~\ref{Basin:Odd3}(b).
Moreover, the dark and light regions are reflections of each other with respect
to the axes of symmetry of $\mathbf T$, where an axis of symmetry of a regular
$n$-gon with odd $n$ is a line connecting a vertex to the midpoint of the
opposite side. In other words, the two basins of attraction are invariant under
the rotational symmetries of the dihedral group $\mathbf D_3$, and are the
images of each other under the reflection symmetries of $\mathbf D_3$. Of
course, the same is true for the attractors $\mathbf T_r$ and $\mathbf T_l$.

\begin{figure}[ht!]
\centering
\subfigure[]{
\includegraphics[width=3cm,height=3cm]{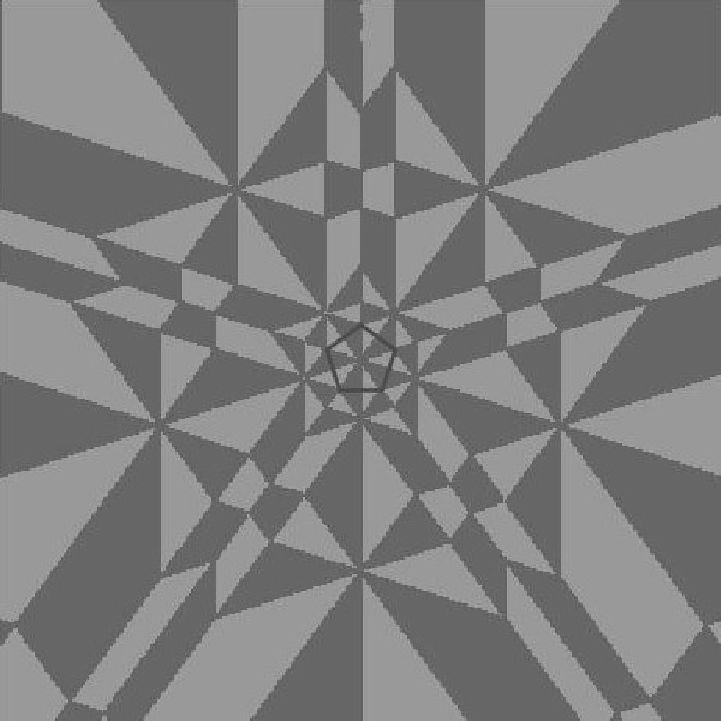} 
}
\hskip 2cm
\subfigure[]{
\includegraphics[width=3.05cm,height=3.05cm]{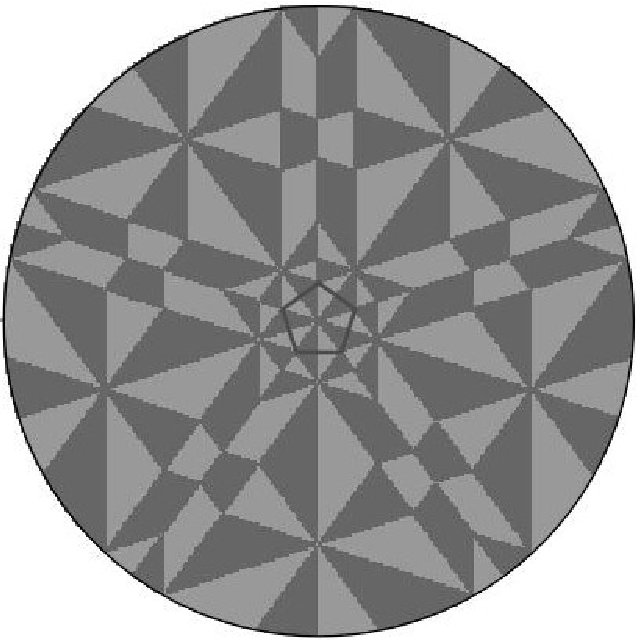} 
}
\caption{Basins for $n=5$ (a) square of side 17 (b) circle of radius
7.} \label{Basin:Odd5}
\end{figure}

If we do the same for a regular pentagon, we obtain the wallpaper shown in
Figure~\ref{Basin:Odd5}. Again, the basins of attraction have the rotational
symmetries of $\mathbf D_5$ and are mirrors of each other with respect to
reflection symmetries of $\mathbf D_5$.
This is particularly interesting as it provides us with quasi pentagonal tiling
of the plane, since a truly pentagonal tiling is forbidden.  Clearly, we can
generate similar tessellations of the plane by finding the basins of attractions
of the FVG played on any regular $n$-gon with odd $n$.

\begin{figure}[ht!]
\centering
\centering
\subfigure[]{
\includegraphics[width=3cm,height=3cm]{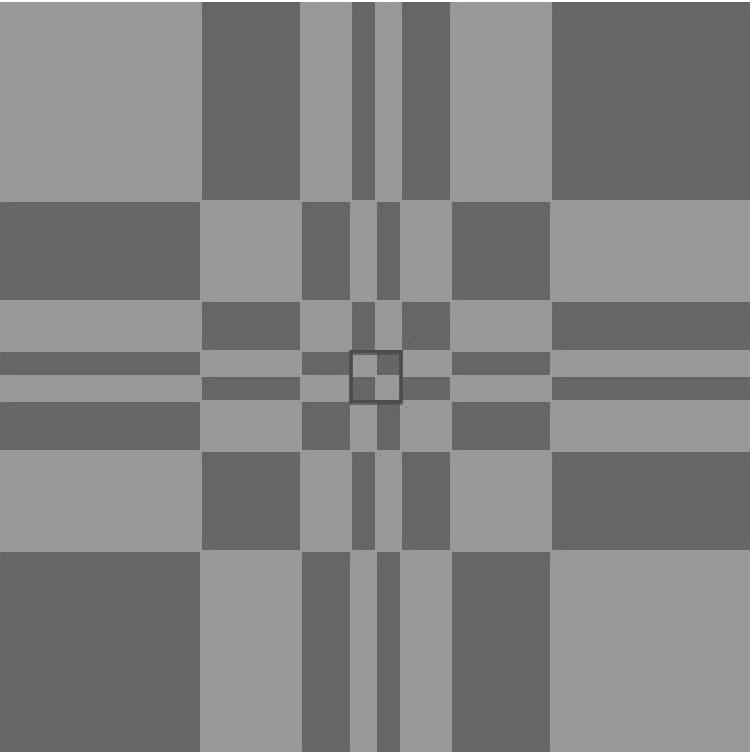}
}
\hskip 2cm
\subfigure[]{
\includegraphics[width=3cm,height=3cm]{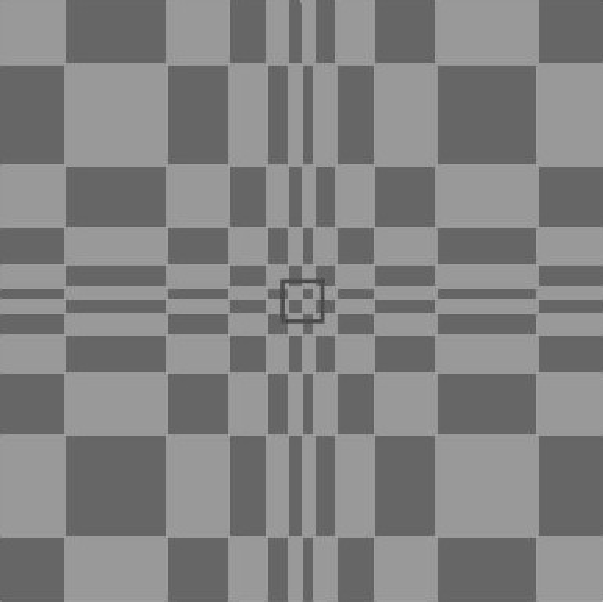}
}
\caption{Basins for $n=4$ in a square of side 15 (a) $r=0.5$ (b)
$r=0.4$.} \label{Basin:Even4}
\end{figure}

For $n$ even, the number of attractors is $n/2$ and so the number of basins
increases with $n$ leading to more intricate wallpapers. Moreover, there are now
$n/2$ axes of symmetry connecting the midpoints of opposite sides and $n/2$ axes
of symmetry connecting opposite vertices.  Taking $n=4$, we get the tiling in
Figure \ref{Basin:Even4}(a). Observe that the basins are now invariant only under a
rotation by $k\pi$, $k=1,2$, and not the full rotational symmetries of the
dihedral group $\mathbf D_4$. This loss of rotational symmetries is compensated
by a gain in reflection symmetries, since in this case each basin is symmetric
with respect to the axes joining opposite vertices.
That is, each basin has both rotational and reflection symmetries. Moreover, the
basins are the images of each other under a rotation by $k\pi/2$, $k=1,3$, and
when reflected with respect to the axes of symmetries obtained by joining the
midpoints of opposite sides.

\begin{figure}[ht!]
\centering
\subfigure[]{
\includegraphics[width=3cm,height=3cm]{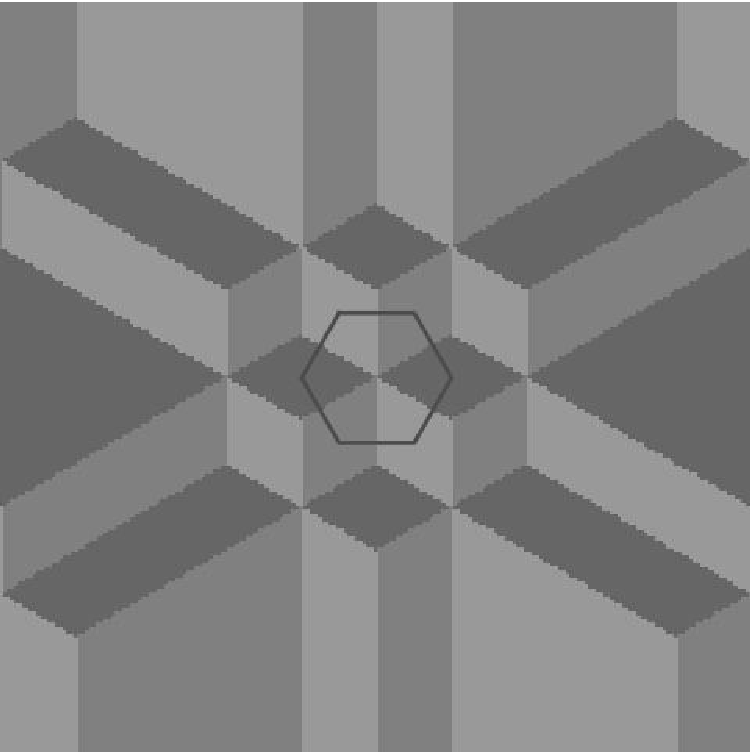}
}
\hskip 2cm
\subfigure[]{
\includegraphics[width=2.7cm,height=3cm]{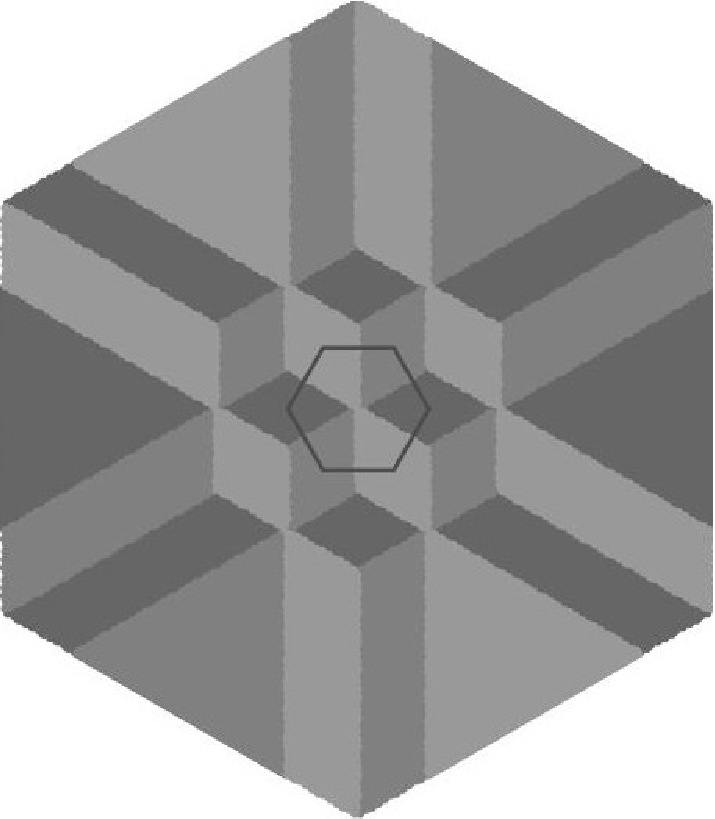}
}
\caption{Basins for $n=5$ (a) Square of side 10 (b) $(x,y)$ with $x\in
[-5,5]$.} \label{Basin:Even6}
\end{figure}

For the regular hexagon $\mathbf P_6$, we get the wallpaper shown in
Figure~\ref{Basin:Even6}(a).  From darkest to lightest, let the three basins of
attraction be called $\mathcal B_1, \ \mathcal B_2$ and $\mathcal B_3$. Also,
let $l_i$ be the axis of symmetry joining the opposite vertices of $\mathbf P_6$
that are contained in $\mathcal B_i$, and let $m_i$ be the axis of symmetry
perpendicular to $l_i$. Clearly, $m_i$ passes through the midpoints of opposite
sides of $\mathbf P_6$ and the symmetry of the tile in
Figure~\ref{Basin:Even6}(b) can be described as follows:
\begin{enumerate}
\item $\mathcal B_i$ is invariant under a rotation by $2k\pi/6$ for $k=3,6$.
\item $\mathcal B_i$ is invariant under a reflection with respect to $l_i$ and
$m_i$.
\item If a symmetry of the dihedral group $\mathbf D_6$ does not preserve
$\mathcal B_i$, then it will map it onto another basin $\mathcal B_j$.
\end{enumerate}

The tile in Figure \ref{Basin:Even6}(b) was generated by taking all points in
$\mathcal B_1$ whose $x$-coordinates are in $[-5,5]$, and then rotating those
points by $\pi/3$ and $2\pi/3$. Of course, the symmetries of the tile, as listed
above, can be extended to the tessellation of the full plane. Moreover, these
symmetries can be similarly described for any regular $n$-gon with $n$ even. In
fact, we can do so without even knowing how the $n/2$ basins look like, since
the symmetries of the basins are directly inherited from the symmetries of the
$n/2$ attractors.

Finally, changing the contraction ratio $r$ amounts to scaling the basins of
attraction, but does not change the symmetries of the tessellation.  In Figure
\ref{Basin:Even4}(b), we changed the contraction ratio from $k_4=0.5$ to
$0.4$. Observe how the overall structure of the wallpaper is preserved, but
everything is now scaled down.


\section{Conclusion} We have seen how random and uniform motion around a regular
$n$-gon produced strikingly different results. Also, we were able to fully
describe the behavior of the URG game and the FVG and in the process identify
their similarities and differences. In particular, we have shown how for $n=3$
the two games lead to the same eventual behavior. That is, they have the same
attractors. But what is most important about this work is that starting with any
initial point in the plane, we can tell exactly where each dynamical system will
end up.

In the URG, there are two attractors for each $n$ and $r$, one for clockwise
rotation and the other for counterclockwise rotation. The attractors can be
determined without iterating the system since they are merely a rotation of the
original $n$-gon by angle $\alpha$ and a scaling down by a factor $\lambda$,
where $\alpha$ and $\lambda$ can be fully calculated using only $n$ and $r$.

For the FVG, the situation is slightly more involved. If $n$ is odd, we still
get two attractors, but we will not know where the initial point will end up
until it enters the original $n$-gon, as the game unfolds. On the other hand, if
$n$ is even, we get $n/2$ attractors, each consisting of a period 2
orbit. Again, we cannot know in advance to which of the attractors will the
initial point converge. To overcome this lack of knowledge about the final
attractor, we iterated a grid of points around the original $n$-gon, which led
us to the discovery of some beautifully intricate patterns that tile the plane.

\bibliography{refer}
\bibliographystyle{amsplain}
\end{document}